\documentclass{article} 
\usepackage{euscript,amsfonts}
\usepackage[final]{epsfig}
\usepackage{psfrag}

\def\nn{\nonumber}

\def\R{\mathbb{R}} 
\def\S{\mathsf{S}}

\def\C{\mathbb{C}}
\def\l{\ell}
\def\hC{\mathsf{C}}
\def\CP{\mathbb{CP}^{n-2}}

\newtheorem{theorem}{Theorem}%[section]
\newtheorem{lemma}[theorem]{Lemma}%[section]
\newtheorem{proposition}[theorem]{Proposition}%[section]
\begin{document}
\newenvironment{proof}{\noindent {\bf Proof.}}{ \hfill$\Box$\\ }
\newenvironment{proofof}[1]{\noindent {\bf Proof of #1.}}{ \hfill$\Box$\\ }

\title{Dual Billiards, Fagnano Orbits, and Regular Polygons.}  
\author{Serge Troubetzkoy}
\date{}

\maketitle

\pagestyle{myheadings}

\markboth{Dual billiards}{SERGE TROUBETZKOY}

%\section{Introduction}
In this article we consider the dual version of two
results on polygonal billiards. We begin by describing these
original results.
The first result is about the dynamics of the so-called pedal map related
to billiards in a triangle $P$. 
The three altitudes of $P$ intersect the opposite
sides (or their extensions) in three points called the feet.  These
three points form the vertices of a new triangle $Q$ called 
the pedal triangle of the triangle $P$ (Figure 1).
It is well known that for acute triangles the pedal triangle forms a
period-three billiard orbit often referred to as the Fagnano orbit,
i.e., the polygon $Q$ is inscribed in $P$ and satisfies
the usual law of geometric optics
(the angle of incidence equals the angle of reflection) or
equivalently (this is a theorem)  the pedal triangle has least perimeter among
all inscribed triangles.
The name Fagnano is used since in 1775 J.\ F.\ F.\ Fagnano gave
the first proof of the variational characterization.
In a sequence of elegant and entertaining articles, J.\ Kingston and
J.\ Synge \cite{KS}, P.\ Lax \cite{L}, P.\ Ungar \cite{U}, and J.\
Alexander \cite{A} studied the dynamics of the pedal map given by 
iterating this process.  The second result, due to DeTemple and
Robertson \cite{DR}, 
is that  a closed convex polygon
$P$ is regular if and only if $P$ contains a periodic billiard path
$Q$ similar to $P$.

There is a dual notion to billiards, called dual or outer billiards.
The game of dual billiards is played outside the billiard
table.  Suppose the table is a  polygon $P$
\begin{figure}[h]\label{fig1}
\centerline{\psfig{file=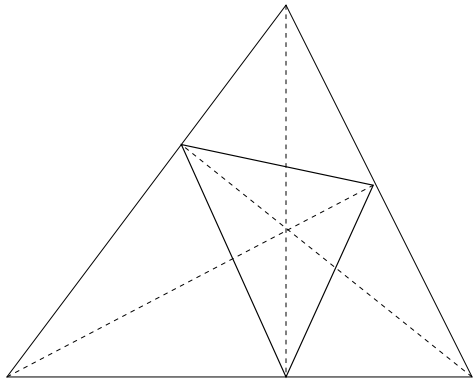,height=50mm}}
\caption{A pedal triangle.}
\end{figure}  
and
that $z$ is a point outside $P$ and not on the continuation of
any of $P$'s sides.   
A line $L$ is a {\em support line} of $P$ if it intersects the boundary
$\partial P$ of $P$ and $P$ lies entirely in one of the two regions
into which the line $L$ divides the plane.
There are two support lines to $P$ through $z$;
choose the right one as viewed from $z$.  If $z$ is not on the
continuation 
of a side of the convex hull of $P$ then this support line intersects $P$ at a single
point which we call the {\em support vertex} of $z$.  Reflect $z$ in
its support vertex to obtain $z$'s image under the dual billiard
map denoted by $T$ (see Figure 2).  
\begin{figure}\label{fig2}
\psfrag*{x}{\footnotesize$z$}
\psfrag*{T(x)}{\footnotesize$Tz$}
\psfrag*{T2(x)}{\footnotesize$\quad\quad T^2z$}
\centerline{\psfig{file=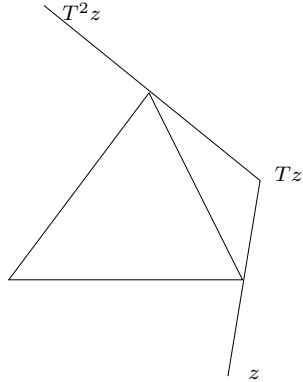,height=50mm}}
\caption{The polygonal dual billiard map.}
\end{figure} The map $T^n$ is defined for all $n$ at the 
point $z$ if none of its images belongs to the continuation of a side of the
convex hull of the polygon.   
%Here we allow self intersecting polygons which we will
%refer to as star shaped polygons. 
Dual billiards have been extensively studied by S.\ Tabachnikov
(see \cite{T}--\cite{TD} and the reference therein).

Throughout the article we will identify  the polygon having vertices 
$z_i \in \C$, $i=1,\dots,n$, ordered cyclically with the 
point $z=(z_1,\dots,z_n)
\in \C^n$.  In particular a polygon for us is an oriented object.
The notion of polygon includes self-intersecting polygons 
(which we call star shaped polygons) and 
geometrically degenerate $n$-gons: those with a side of length 0, or an angle of $\pi$
or $2 \pi$ at a vertex.   Corresponding to the fact that the Fagnano
billiard orbit hits each side of the triangle, we call 
an $n$-periodic dual billiard orbit consisting of $n$ points consecutively reflected
in the vertices $z_1,z_2,\dots,z_n$ a
{\em Fagnano} dual billiard orbit. 
Fagnano orbits for dual billiards were introduce by Tabachnikov in
\cite{T1}, where he studied a certain variational property analogous
to one in the pedal case studied by Gutkin \cite{G}.

Motivated by the notion of pedal triangles, we introduce here dedal $n$-gons.
An $n$-gon $Q$ is called a {\em dedal} $n$-gon of the $n$-gon
$P$ if reflecting the vertex $w_i$ of $Q$ in
the vertex  $z_i$ of $P$ yields the vertex
$w_{i+1}$. (Throughout the article all subscripts will be
taken modulo $n$ without explicit mention.)
There is no requirement that the sides of $Q$ touch
$P$ only at a vertex.
From the definition it is clear that if
a Fagnano orbit of an $n$-gon $P$ exists then it is a dedal $n$-gon. 
A nondegenerate polygon can have a degenerate dedal
polygon and vice-versa.  Examples are shown in Figures 3 and
4.  The degeneracy not shown cannot occur: i.e.,
it is impossible
for two consecutive vertices of $P$ to coincide if $Q$ is
nondegenerate. We will not dwell on this aspect.

%The space of $n$-gons $Q$ with degenerate $\mu(Q)$ is
%a complex co-dimension 1 subspace.
\begin{figure}
\psfrag*{w1}{\footnotesize$w_1$}
\psfrag*{w2}{\footnotesize$w_2$}
\psfrag*{w3}{\footnotesize$w_3$}
\psfrag*{w4}{\footnotesize$w_4$}
\psfrag*{w5}{\footnotesize$w_5$}
\psfrag*{z1}{\footnotesize$z_1$}
\psfrag*{z2}{\footnotesize$z_2$}
\psfrag*{z3}{\footnotesize$z_3$}
\psfrag*{z4}{\footnotesize$z_4$}
\psfrag*{z5}{\footnotesize$z_5$}
\psfrag*{a}{\footnotesize$z_2=z_5$}
\centerline{\psfig{file=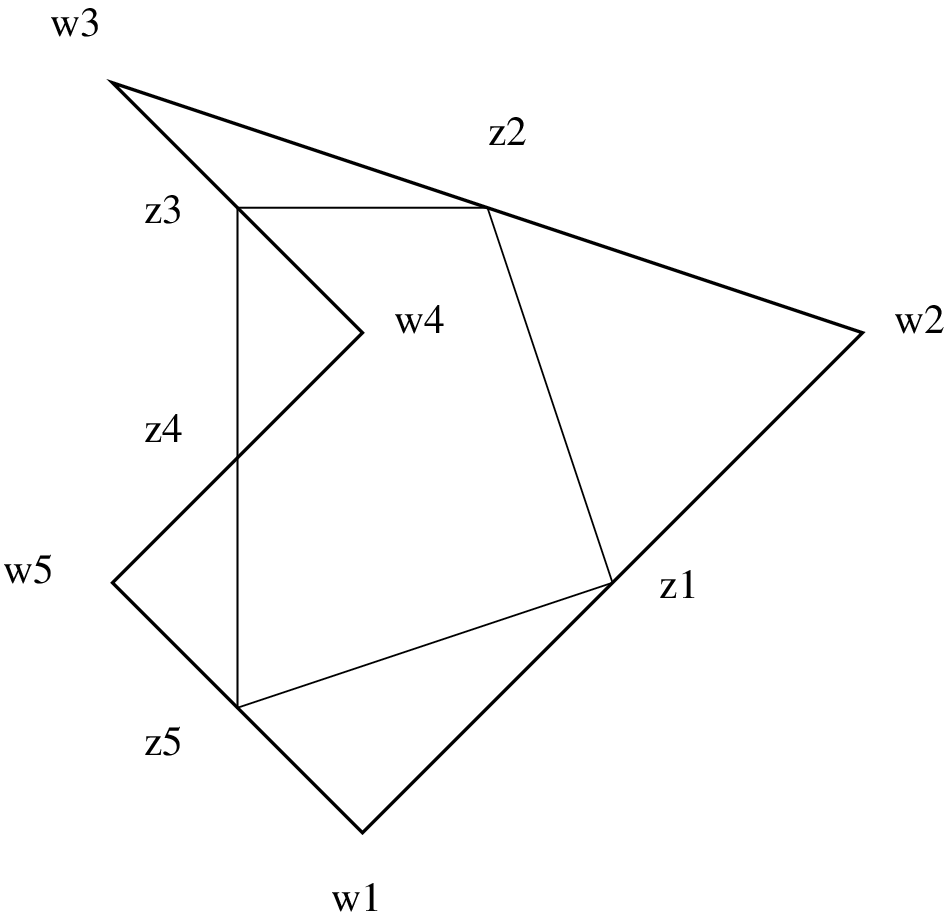,height=50mm} \hspace{0.5cm}
\psfig{file=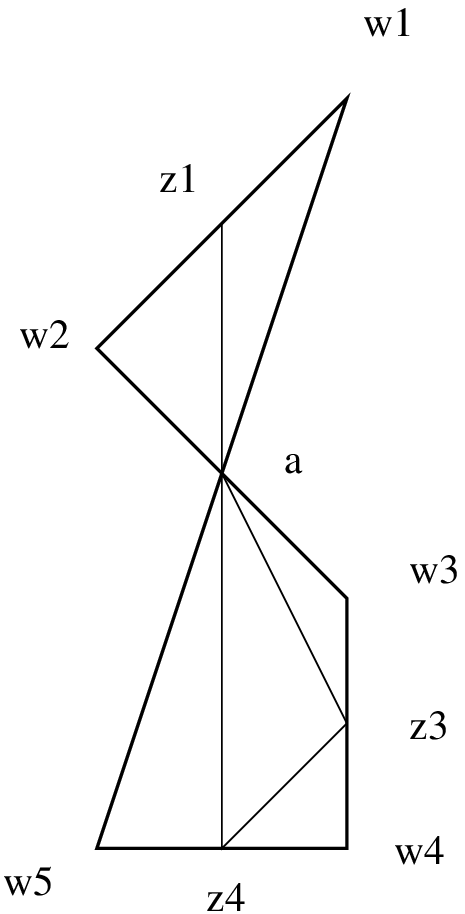,height=50mm}}
\caption{Nondegenerate dedal pentagons corresponding to degenerate 
pentagons, with  angle $\pi$ and $2\pi$.}
\end{figure}\label{fig6}

\begin{figure}\label{f5}
\psfrag*{w1}{\footnotesize$w_1$}
\psfrag*{w2}{\footnotesize$w_2$}
\psfrag*{w3}{\footnotesize$w_3$}
\psfrag*{w4}{\footnotesize$w_4$}
\psfrag*{w5}{\footnotesize$w_5$}
\psfrag*{z1}{\footnotesize$z_1$}
\psfrag*{z2}{\footnotesize$z_2$}
\psfrag*{z3}{\footnotesize$z_3$}
\psfrag*{z4}{\footnotesize$z_4$}
\psfrag*{z5}{\footnotesize$z_5$}
\psfrag*{aa}{\tiny$w_3=w_4$}
\psfrag*{a}{\tiny$=z_3$}
\centerline{\psfig{file=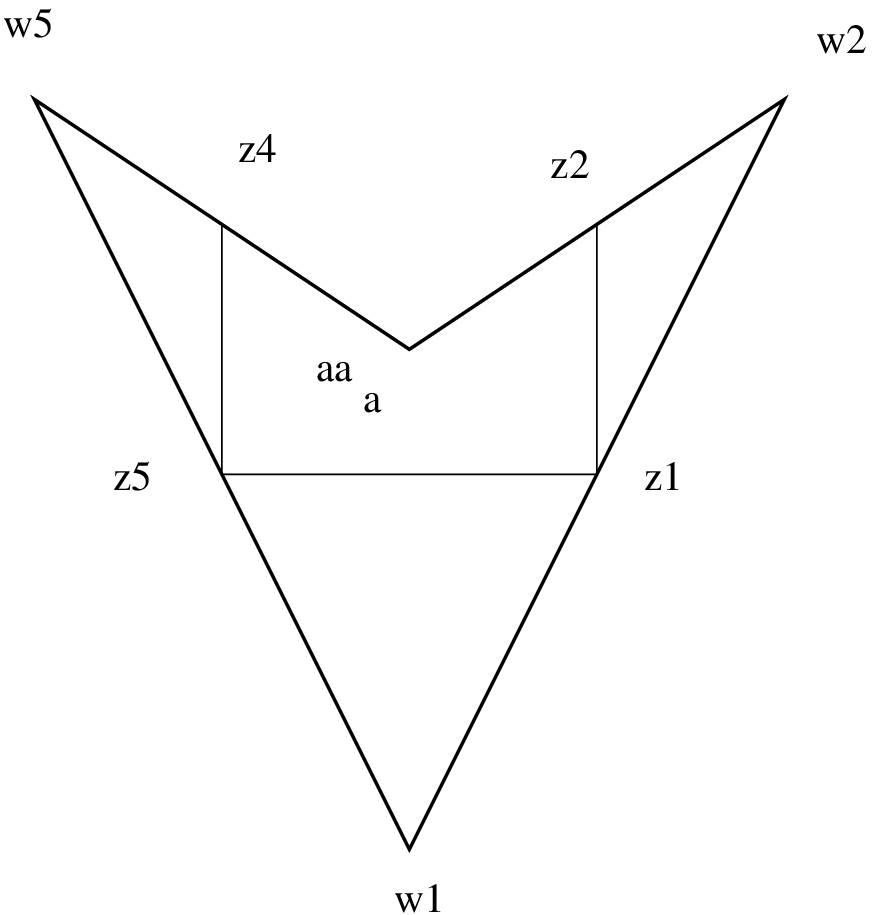,height=34mm} \hspace{0.4cm}
\psfig{file=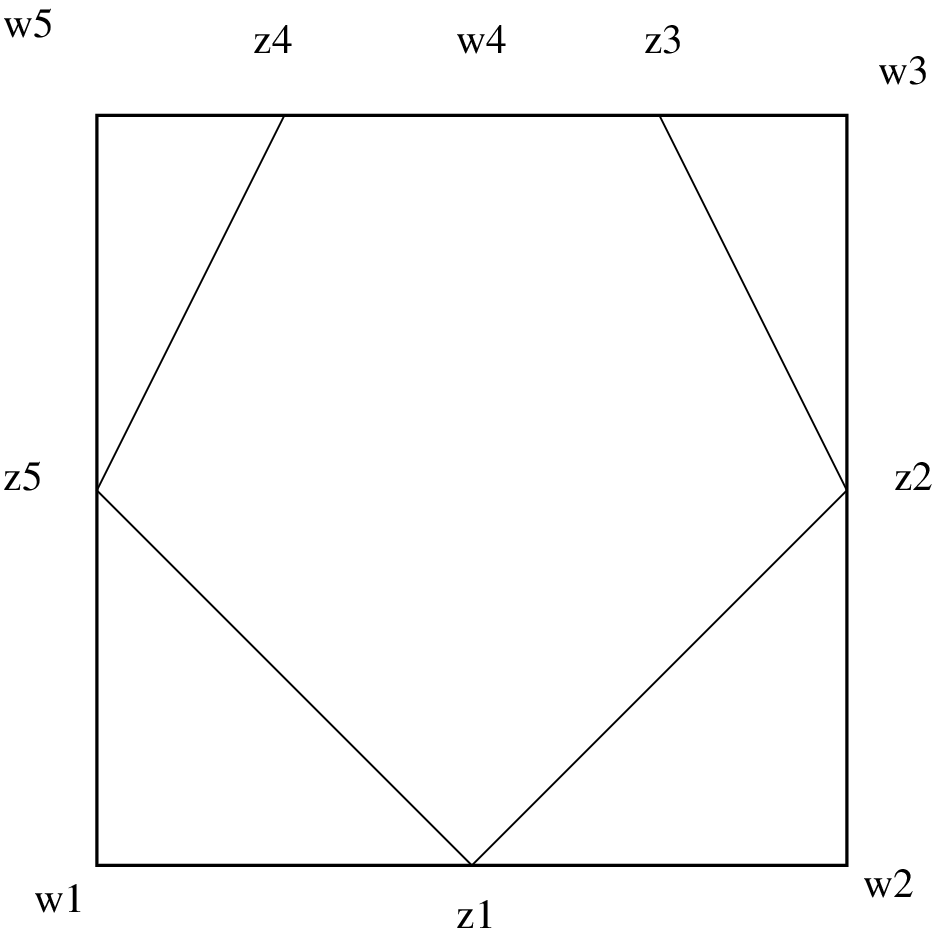,height=38mm} \hspace{0.4cm}
\psfig{file=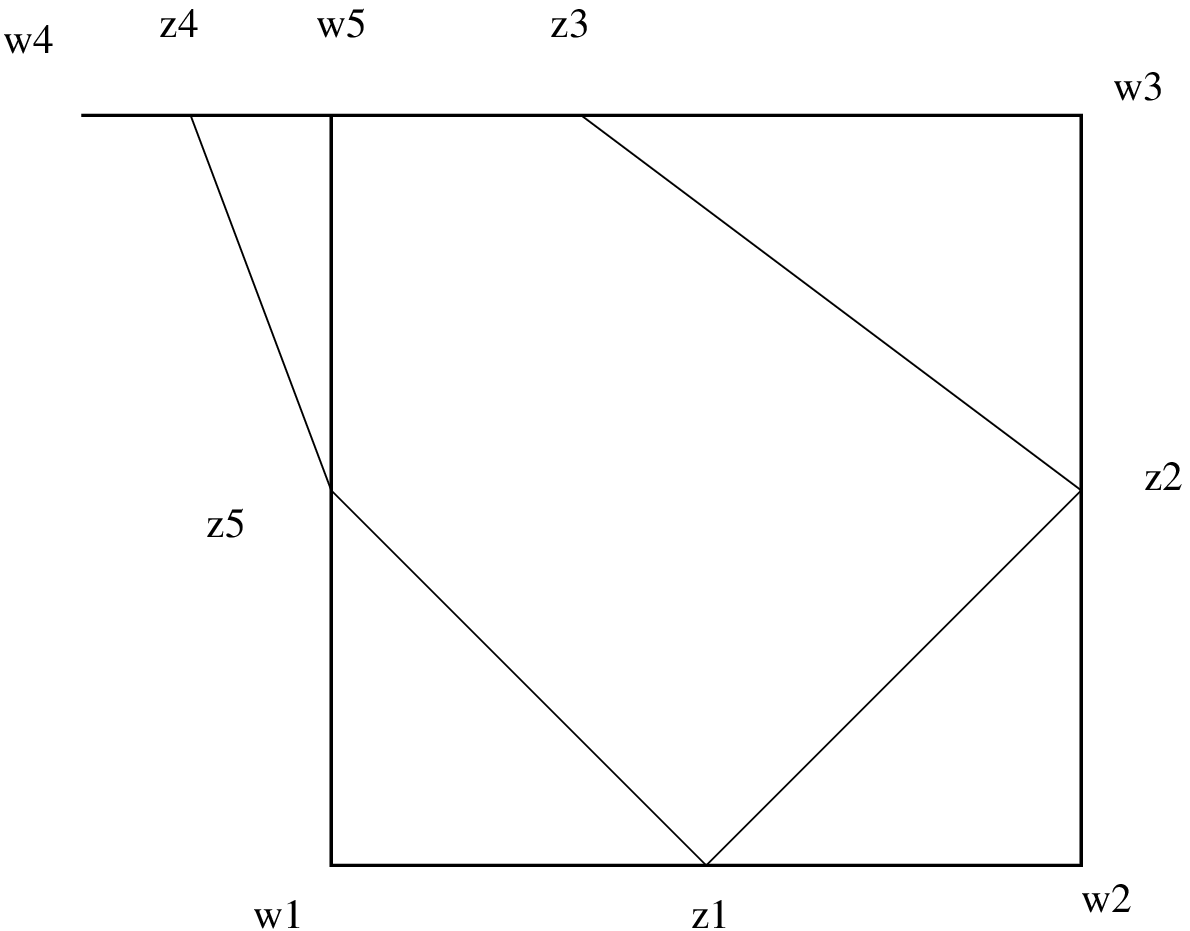,height=38mm}}
\caption{Degenerate dedal polygons: loss of a vertex, angle $\pi$,
and  angle $2\pi$.}
\end{figure}

In this article we study the dedal $n$-gon(s) $Q$ of an $n$-gon $P$.
Our main results are the following.   If $n$ is odd then its dedal 
$n$-gon exists and is unique.  For $n$ even we give a necessary
and sufficient condition for the existence of dedal $n$-gons and
describe the space of all dedal $n$-gons of $P$.
Then we go on to characterize regular and affinely regular $n$-gons
by similarity to their dedal $n$-gons.
Finally we give a complete description of the dynamics of the
{\em developing map} $\mu(Q) :=  P$.
The proofs of all our results boil down to some
linear algebra of the dedal map.

After we wrote this article one of the anonymous referees pointed
out that iteration of the developing map has already been studied by
Berlekamp, Gilbert, and Sinden in 1965 \cite{BGS}.  They
answer a question they attribute to G.\ R.\ MacLane,
namely, they prove that for almost every  polygon $Q$
there exists an $M \ge 1$ such that $\mu^M(Q)$ is convex. Note that the
image of a convex polygon is convex; thus this implies that for almost every
$Q$ there exists an $M$ such that $\mu^m(Q)$ is convex for all $m \ge
M$. 

\newpage

Suppose $Q(w_1,\dots,w_n)$ is a dedal polygon of $P(z_1,\dots,z_n)$.  
The definition implies that 
$$
z_i = (w_i + w_{i+1})/2 
$$
(see Figures 3, 4, and 5).
The linear transformation $\mu: \C^n \to \C^n$  given by
$\mu(w_1,\dots,w_n) = (z_1,\dots,z_n)$ is called the  
{\em developing map}.

The characteristic polynomial of $\mu$ is $$(1- 2x)^n - (-1)^n.$$
Its eigenvalues are $(1 + q^i)/2$ for $i=0,1,\dots,n-1$, where
$q := \exp(2\pi i/n)$.
All the eigenvalues differ from zero except for the $(n/2)$th  eigenvalue
when $n$ is even. The $i$th eigenvector is  
$X_i := (1,q^i,q^{2i},\dots,q^{(n-1)i})$.
%A simple calculation (see \cite{BGS},\cite{T1}) shows that 
The vectors $X_i$ form a basis of our space of polygons.
If $i$ divides $n$ then $X_i$ is a polygon with $n/i$
sides. Nonetheless 
for the sake of clarity (avoiding stating special cases) we will
think of this as an $n$-gon which is traced $i$ times; for example if
$n=6$ then $X_2 = (1,q^2,q^4,q^6,q^8,q^{10}) = (1,q^2,q^4,1,q^2,q^4)$
traces the triangle $(1,q^2,q^4)$ twice.
An exception to this rule is the case $n$ even and $i= n/2$.  In this
case $X_{n/2}$ is a segment which we do not consider as a 
polygon.

The following proposition clarifies
the existence of dedal polygons.
\begin{proposition}\label{prop1}
a) Suppose that $n \ge 3$ is odd.  Then for any $n$-gon $P$ there is a
unique dedal $n$-gon $Q$.\\
b)  If $n \ge 3$ is even then dedal $n$-gons exist if and only if
the vertices of $P$  satisfy $z_1 - z_2 + z_3 -
\cdots - z_n = 0$. If this equation is satisfied then 
there is a unique dedal n-gon $Q_0$ in the space $X_{n/2}^\perp
:= \{z \in \C^n:\   z \cdot X_{n/2} = 0\}$.
The set $D:=\{Q_0 + s X_{n/2}: \ s \in \mathbb{C}\}$ consists of the
dedal $n$-gons of $P$. In particular
for each $i \in \{1,\dots,n\}$ 
every point $w \in \mathbb{C}$ is the $i$th vertex of a unique dedal
$n$-gon  $Q_i(w)$.
\end{proposition}
\noindent
We remark that the condition $z_1-z_2+ z_3 - \cdots - z_n = 0$ means that
the center of mass of the even vertices coincides with the center of
mass of the odd vertices.
  
\begin{proofof}{Proposition \ref{prop1}}
If $n \ge 3$ is odd then the map $\mu$ is invertible
with 
$$
w_i = z_i - z_{i+1} + z_{i+2} - \cdots + z_{i-1},
$$
and thus the dedal polygon exists and is unique.
%We can also see the uniqueness geometrically, the composition of an
%odd number of reflections is again a
%reflection, thus a dedal polygon is isolated and unique.

On the other hand, if $n \ge 3$ is even then the map $\mu$
is not invertible.  The kernel
is one (complex) dimensional and is generated by the vector 
$X_{n/2} := (1,-1,1,-1,\dots,1,-1)$.  Dedal
polygons  exist if and only if 
$z = (z_1,z_2,\dots,z_n)$ is in the range of $\mu$, i.e.,
the space spanned by the $X_i$
for $i \ne n/2$. The range and kernel of $\mu$ are orthogonal
since if $i \ne n/2$ then
$$X_i \cdot X_{n/2} = (1 + q^{2i} + q^{4i} + \cdots + q^{(n-2)i}) 
- (q^i + q^{3i} + \cdots + q^{(n-1)i}) = 0-0 = 0.$$
Thus $z$ is in the range of $\mu$ if and only if
it satisfies $z \cdot X_{n/2} = 0,$ or equivalently  
\begin{equation}\label{existence}
z_1 - z_2 + z_3 - \cdots - z_n = 0.
\end{equation}
Alternatively, to see this note that 
\begin{eqnarray}
&\hspace{-10pt}z_1 + z_3 + \cdots + z_{n-1} 
= \frac12 \left (   (w_1 + w_2) + (w_3 + w_4) + \cdots + (w_{n-1} +
  w_n) \right )\nonumber\\
& = \frac12 \left ( (w_2 + w_3) + (w_4 + w_5) + \cdots + (w_n + w_1) \right ) 
= z_2 + z_4 + \cdots + z_n.\nonumber
\end{eqnarray}
The uniqueness of $Q_0$ follows since the map $\mu$ is invertible 
on the space $X_{n/2}^\perp$.  
The statement about the set $D$ follows immediately since
$X_{n/2}$ is the kernel of $\mu$. Let $Q_0 := (w_1^0,\dots,w_n^0)$.
For each $i$, any point $w \in \mathbb{C}$
can be uniquely expressed as $w_i^0 + s$ for some $s \in \mathbb{C}$.
\end{proofof}
 
Suppose that a polygon $Q$ without self intersection is a dedal
polygon of $P$.
If $Q$ is convex then
it is clearly a Fagnano orbit of $P$ since by convexity the
polygon $P$ must be contained in $Q$. 
On the other hand if $Q$ is not
convex then it cannot be a Fagnano orbit since $P$ cannot be contained
in $Q$.
In particular Fagnano orbits 
always exist for triangles, but not for polygons with more sides.
An example of a pentagon without a Fagnano orbit is given
in Figure 5. Although not every polygon  has a Fagnano orbit, 
it does have a periodic orbit that is  nearly as simple. Namely, consider the
second iteration $T^2$ of the dual billiard map. Connecting  
the consecutive points of a periodic orbit of $T^2$
yields a polygon.  Cutler has shown that the map $T^2$ has a periodic
orbit
which lies outside any compact neighborhood of $P$ and the polygon
constructed from the orbit makes a
single turn about $P$ \cite{T2}. The existence of a periodic
orbit for the usual billiard in an arbitrary polygon remains open.
 
\begin{figure}
\psfrag*{w1}{\footnotesize$w_1$}
\psfrag*{w2}{\footnotesize$w_2$}
\psfrag*{w3}{\footnotesize$w_3$}
\psfrag*{w4}{\footnotesize$w_4$}
\psfrag*{w5}{\footnotesize$w_5$}
\psfrag*{z1}{\footnotesize$z_1$}
\psfrag*{z2}{\footnotesize$z_2$}
\psfrag*{z3}{\footnotesize$z_3$}
\psfrag*{z4}{\footnotesize$z_4$}
\psfrag*{z5}{\footnotesize$z_5$}
\centerline{\psfig{file=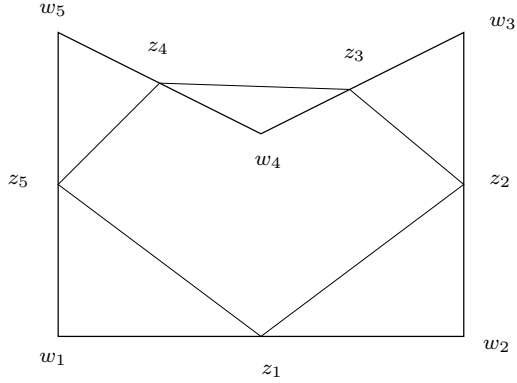,height=50mm}}
\caption{ A pentagon without a Fagnano orbit. The unique dedal pentagon
is pictured.}\label{fig3}
\end{figure} 

The set of polygons with
center of mass at the origin is
$$\hC := \{(w_1,\dots,w_n) \in \C^n: \ w_1 + \cdots + w_n =~0\}.$$
The eigenvector $X_{0}  = (1,1,\dots,1)$ represents a   
polygon which degenerates to a
point.  The set $\hC$ is the orthogonal complement of $X_0$, i.e.,
$\hC = \{w=(w_1,\dots,w_n) \in \C^n: \  w \cdot X_0 = 0\}$.
We will express polygons in the eigenbasis,
i.e., we write $P = \sum_{i=0}^{n-1} a_i X_i$; the coefficients $a_i$ are
complex numbers. Since the map $\mu$ preserves the center of mass,
throughout the rest of the article we assume that the center of mass is
at the origin, i.e., $a_0=0$.

$X_1$ and $X_{n-1}$ are the
usual regular $n$-gon in counterclockwise and clockwise orientation
(Figure 6a).
If $i$ and $n$ are relatively prime and $i \not \in \{1,n-1\}$ then
$X_i$ is star shaped and we also call it regular (Figures 6b and c).
Finally if $i$ divides $n$ and $i \ne n/2$ then $X_i$ is naturally
a regular $n/i$-gon which, as mentioned above, 
we will think of as (a multiple cover of)
a regular $n$-gon.
\begin{figure}\label{fig4}
\psfrag*{a}{\footnotesize$(a)$}
\psfrag*{b}{\footnotesize$(b)$}
\psfrag*{c}{\footnotesize$(c)$}
\centerline{\psfig{file=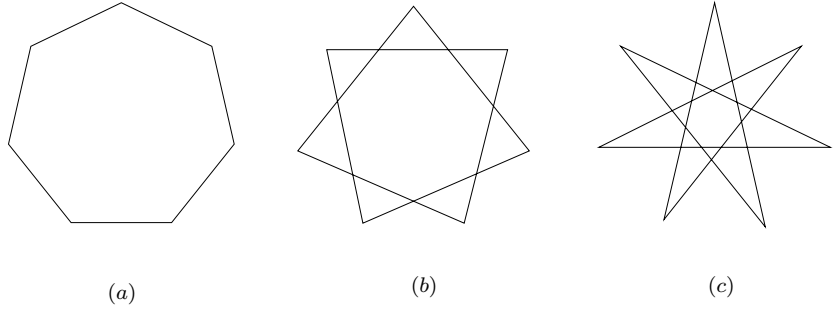,height=40mm}}
\caption{Up to orientation there are three regular 7-gons.}
\end{figure}  

Two (unoriented) polygons are called {\em similar} if all corresponding angles
are equal and all distances are increased (or decreased) in the same
ratio. 
%Since we have assumed that the center of mass is at the origin
%two polygons $P$ and $Q$ are similar iff there is a nonzero complex 
%constant $k$ such that 
%$P = k \, Q$. 
Since we study oriented ordered polygons we additionally want
the marking of the vertices of $P$ and $Q$ to correspond, in which case
we will call $P$ and $Q$ {\em $\star$-similar}. Thus  
two polygons
$P = \sum_{i=1}^{n-1} b_i X_i$ and 
$Q = \sum_{i=1}^{n-1}  a_i X_i$  are
$\star$-similar if there exists a nonzero complex
constant $\l$ such that $b_i = \l \, a_i$ for all $i$.
We will also write this as $P = \l \, Q$.
Note that if $P$ and $Q$ are $\star$-similar then they are similar.
On the other hand if $P$ and $Q$ are similar then $P$ is $\star$-similar
to a cyclic permutation 
$Q^{(k)} := (w_k, w_{k+1},w_{k+2}, \dots, w_{k-1})$ of $Q$ or
a cyclic permutation 
$\bar{Q}^{(k)} := (w_k, w_{k-1}, w_{k-2} ,\dots, w_{k+1})$
of $Q$ with the opposite orientation.

In analogy to DeTemple
and Robertson's result we have:
\begin{theorem}\label{thm0}
Fix $n \ge 3$.  An $n$-gon $P$ is regular if and only if it has
a dedal polygon $Q$ which is $\star$-similar to $P$.
\end{theorem}
Note that if $n$ is odd then 
$Q$ is the unique dedal polygon
of $P$, while if $n$ is even then $Q$ is the unique dedal polygon
$Q_0 \in X_{n/2}^\perp$.

\begin{proof}
Suppose $P$ is regular, i.e.,  there is a nonzero complex constant $\l$
such that $P = \l \, X_j$, where $j \ne n/2$
if $n$ is even. 
If $n$ is odd then let $Q$ be the unique dedal polygon of $P$. 
If $n$ is even then since $P$ is regular it satisfies
(\ref{existence}). Thus $P$ has dedal polygons and
we choose $Q = Q_0 \in X_{n/2}^\perp$.
In both cases let $Q=\sum a_i X_i$. Since 
$P=\mu(Q) = \sum \frac{1+q^i}{2} a_i X_i$,
we have  $a_i = 0$ for $i \ne j$, i.e., $Q = a_j  X_j$. Thus $Q$ is 
$\star$-similar to $P = \ell \, X_j$.

Conversely suppose that $Q = \sum a_i X_i$ is a dedal
polygon of $P = \sum b_i X_i$ and that $P$ and $Q$  
are $\star$-similar, i.e., there is a nonzero complex constant $\l$ such
that $a_i = \l \, b_i$. 
%Note that since  $Q$ exists either $n$ is odd or that $n$ is even and
%$P$ satisfies \eqref{existence}.
Since $P = \mu(Q)$ we have $\frac{1+q^i}{2} a_i =  b_i$ for $i=1,\dots,n$.
Combining this with the $\star$-similarity yields
$\frac{1+q^i}{2} = 1/{\l}$ (a constant) for each $i$ such that $a_i \ne 0$.  
It follows $q^i = q^j$ if $a_i \ne 0$ and $a_j \ne
0$. Since $q^i \ne q^j$ for $i \ne j$ it follows that  only a single
$a_j$ is nonnull. Therefore $Q = a_j X_j$ and thus $P = \ell^{-1} \,
a_j X_j$, and both  are regular. 
Note that if $n$ is even then $i \ne n/2$, since
if $i = n/2$ then  $(1 + q^i)/2 = 0$ and thus $b_i = a_i = 0$.
\end{proof}

Of course we would like to have the analog of Theorem \ref{thm0} with 
the usual notion of similarity. We call an $n$-gon {\em  affinely regular} 
if there exists a $j$ ($j \ne n/2$ when $n$ is even)  
such that $P$ is in the subspace $A_j$ generated by $X_j$ and $X_{n-j}$.
This is equivalent to $P$ being the image of a regular $n$-gon by an
affine transformation, which explains the name.

 \newcounter{Lcount}

\begin{theorem}\label{thm1} a) Suppose $n \ge 3$ is odd.
An $n$-gon $P$ is affinely regular if and only if 
it has a dedal polygon $Q$ which is similar to $P$.\\
b) If $n \ge 4$ is even then an $n$-gon $P$  appears in the following
list of affinely regular polygons if and only if
it has a dedal polygon $Q$ which is similar to $P$.  
\begin{list}{\roman{Lcount})}
  %    inform the list command to use this counter
    {\usecounter{Lcount}
  %    set rightmargin equal to leftmargin
    \setlength{\rightmargin}{\leftmargin}}
  %    we can now begin the "items"
\item Regular $n$-gons 
\item Any $P \in A_j$ such that there exists $k \in
\{1,,\dots,n\} \setminus~\{n/2\}$  such that
$n$ divides $j(2k-1)$
\item Any $P = b_j X_j + b_{n-j}X_{n-j} \in A_j$ such that
there exists $k \in \{1,\dots,n\}$ with
${b_j}/{b_{n-j}} = \pm  q^{j(k+3/2)}$
\end{list}
\end{theorem}
All triangles are affinely regular. Berlekamp et al.\ 
noticed that every triangle is similar to its dedal triangle
\cite{BGS}. One would like to know if there are other
special properties of the polygons in the list.

\begin{proof}
Suppose $Q^{(1)} = Q = \sum_{i=1}^{n-1} a_i X_i$.  Then 
it is a simple exercise in linear algebra to verify that 
$Q^{(k)} = \sum a_i q^{i(k-1)} X_i$ and   
%$ (a_1q^k,a_2q^{2k},\dots,a_{n-1}q^{(n-1)k})$ and 
that $\bar{Q}^{(k)} = \sum a_{n-i} q^{i(k+1)} X_i$.
%$(a_{n-1}q^{k},a_{n-2}q^{2k},\dots,a_1q^{(n-1)k})$.
Thus the similarity of $P$ to $Q$ 
is equivalent to the existence of a $k \in \{1,\dots,n\}$ and
a nonzero complex constant $\l$ such that
$b_i = \l \, a_iq^{i(k-1)}$ for $i = 1,\dots,n-1$
or
$b_i = \l \, a_{n-i}q^{i(k+1)}$ for $i = 1,\dots,n-1$.

The structure of the proof is as follows.  We group the even and odd
cases together and start by proving that
the various classes of affinely regular polygons are similar to
their dedal polygons.  Then we turn to the converse.
 
Suppose that $n$ is odd and
$P = b_j X_j + b_{n-j} X_{n-j}$ is affinely regular.
Let $Q = Q^{(1)} = \sum a_i X_i$
be the unique dedal polygon of $P$.  
Since $P=\mu(Q) = \sum \frac{1+q^i}{2} a_i X_i$,
we have  $a_i = 0$ for $i \not \in  \{j,n-j\}$. 
We claim that there is a nonzero complex constant $\ell$ such that
$P = \ell \, Q^{((n+3)/2)}$.
To see this note that
$$
1 + q^j = q^j(1 + q^{-j}) = q^{j(n+1)} ( 1 + q^{-j}) = q^{j(n+1)/2}
q^{j(n+1)/2} ( 1 + q^{-j})$$
or
$$\frac{1+q^j}{2q^{j(n+1)/2}} = \frac{1+q^{-j}}{2q^{-j(n+1)/2}}.$$
Thus $P = \l  \, Q^{((n+3)/2)}$  for $\l := (1+q^j)/ (2q^{j(n+1)/2})$.

The case  $n$ is even and $P$ of class (ii) is similar.
Consider the dedal polygon $Q = Q_0 \in X_{n/2}^{\perp}$.
Again $P = \mu(Q)$ implies that $a_i = 0$ if $i \not \in \{j,n-j\}$.
We claim that there is a nonzero complex constant $\ell$ such that
$P = \ell \, Q^{(k+1)}$.  Since $n$ divides $j(2k-1)$ we have
$$
1 + q^j = q^j(1 + q^{-j}) = q^{j + j(2k-1)} ( 1 + q^{-j}) = q^{kj}
q^{kj} ( 1 + q^{-j})$$
or
$$\frac{1+q^j}{2q^{kj}} = \frac{1+q^{-j}}{2q^{-kj}}.$$
Thus $P = \l \, Q^{(k+1)}$ for $\l := (1+q^j)/(2q^{kj})$.

The case $n$ even and $P$ regular has already been treated in Theorem
\ref{thm0}; therefore it remains to treat case (iii) of $n$ even.  
We suppose $P$ is of this form and want to show 
that $P$ is similar to $\bar{Q}^{(k)}$, i.e., that there is a nonzero
complex constant $\ell$ such that
$b_i = \ell \, a_{n-i} q^{i(k+1)}$ for all $i$.
As before $P = \mu(Q)$ implies that $a_i = 0$ if $i \not \in
\{j,n-j\}$.
Combining the two relations,
$b_i = \frac{1+q^i}{2} a_i$
and
$b_j/b_{n-j} = \pm  q^{j(k+3/2)}$
yields
$$ b_{j} = \pm q^{j(k+3/2)} b_{n-j}
= \pm q^{j(k+3/2)} \frac{1+q^{-j}}{2} a_{n-j}
= \pm q^{j(k+1)}  \frac{q^{j/2}+q^{-j/2}}{2} a_{n-j}$$
and 
$$ b_{n-j} = \pm q^{-j(k+3/2)} b_j
=  \pm q^{-j(k+3/2)} \frac{1+q^j}{2} a_j
= \pm q^{-j(k+1)} \frac{q^{j/2} +q^{-j/2}}{2} a_j.$$
Choosing $\ell := \pm  (q^{j/2} + q^{-j/2}) /2 \in \R$ yields
$P = \ell \, \bar{Q}^{(k)}$.

We turn now to the converse.  We will first
prove that for $n$ even or odd if $Q = 
\sum a_i X_i$ is similar to $P = \sum b_i X_i$
then $P$ is affinely regular.

We first treat that case when 
$P = \l \, Q^{(k+1)}$ for some $k$ and $\l \in \C \setminus \{0\}$, 
i.e., $b_i = \l \, a_i  q^{ik}$. 
Since $P = \mu(Q)$ we have $\frac{1+q^i}{2} a_i =  b_i$ for $i=1,\dots,n-1$.
Note that if $n$ is even then 
$(1 + q^{n/2})/2 = 0$, and thus $b_{n/2} = \frac{1+q^{n/2}}{2} a_{n/2}
=  0$ and $a_{n/2} = \l^{-1} q^{-nk/2} b_{n/2}$ is zero as well.
For each $i$ such that $a_i \ne 0$, combining the two relations 
between $a_i$ and $b_i$ yields 
$q^{-ik} \frac{1+q^i}{2} = \l$.
If only a single $a_j$ is nonnull, then $Q$ is regular, and since 
$P$ is similar to $Q$ it is regular as well.  
Now suppose that $a_i$ and $a_j$ are nonnull. Then the previous equation
implies that 
\begin{equation}(1+q^i)/(1+q^j) = q^{k(i-j)}.\label{1}\end{equation}
Taking absolute values yields $|(1+q^i)|=|(1+q^j)|$, which implies
$i = \pm j$.  Thus $P$ is affinely regular.

If $n$ is even we need to conclude more.
Note that $(1+q^{-j})/(1+q^j) = q^{-j}$. Thus
taking $i = -j$ in (\ref{1}) implies $1 = q^{j(1-2k)}$.
Thus $j(2k-1)$ is a multiple of $n$, i.e., we are in case (ii) of the
list.
 
Finally suppose that $P= \mu(Q)$ and $Q$ are similar but
have the opposite orientation, i.e., $P = \l \, \bar{Q}^{(k)}$
for some $\l \in \C \setminus \{0\}$, or equivalently 
$b_i = \l \, a_{n-i} q^{i(k+1)}$.
Combing this with $\frac{1+q^i}{2} a_i = b_i$ yields
$\frac{1+q^{i}}{2} a_i = \l \, a_{n-i} q^{i(k+1)}$. 
If $n$ is even and $i = n/2$ then this equation implies that
$a_{n/2} = 0$.  For all other cases 
it implies that
$a_i$ and $a_{n-i}$ are simultaneous zero or nonzero.
For any $i$ 
such that they are nonzero there are two such equations, and they
imply
\begin{equation}\label{2}
\frac{a_i}{a_{n-i}} = \frac{2 \, \l \, q^{i(k+1)}}{1+q^{i}} =   
\frac{1+ q^{-i}}{2 \, \l \, q^{-i(k+1)}} .
\end{equation}
This in turn yields
\begin{equation}\label{3}
4 \, \l^2  = (1+q^i) \cdot (1 + q^{n-i}) = 
{2 + q^i + q^{n-i}} := f(i) .
\end{equation}
It is easy to see that $f(i) = f(j)$ if and only if $j=i$ or $j=n-i$.
Since the left-hand side of Equation (\ref{3}) does not depend on $i$, 
there is exactly one pair $(a_j,a_{n-j})$ of nonzero coefficients,
i.e., $Q$ is affinely regular.  Since $P$ is similar to $Q$ it is also
affinely regular.

Finally suppose that $n$ is even.  Equations (\ref{2}) and (\ref{3}) 
imply that for
each $j$ there exists $k \in \{1,2 \dots ,n\}$ such that
$$\frac{b_j}{b_{n-j}} = \pm  \frac{\sqrt{f(j)} q^{j(k+1)}}{1+q^{-j}}
= \pm \sqrt{\frac{1+q^j}{1+q^{-j}}} q^{j(k+1)}
= \pm q^{j(k+3/2)},$$
i.e., we are in case (iii) of the list.
\end{proof}
 
\begin{lemma}\label{n}
If $Q$ is affinely regular then $\mu^{n}(Q)$ is $\star$-similar to $Q$.
\end{lemma}

\begin{proof}
Suppose $Q = a_j X_j  + a_{n-j}X_{n-j}$.
Then
$$\mu^{n}(Q) = \Big (\frac{1+q^j}{2} \Big )^{n} a_j X_j +  \Big
(\frac{1+q^{n-j}}{2} \Big )^{n} a_{n-j} X_{n-j}.$$
But
\begin{eqnarray}\label{aha}
(1+q^{n-j})^n &=& \sum_{k=0}^n {n  \choose k} q^{(n-j)k} = 
\sum_{k=0}^n {n  \choose n-k} q^{-jk}  \\ 
&=& \sum_{i=0}^n {n  \choose  i} q^{-j(n-i)} 
= \sum_{i=0}^n {n \choose  i}  q^{ji} =(1+q^j)^n.\nn
\end{eqnarray}
Thus $\mu^n(Q) = \l \, Q$ with $\l = (\frac{1+q^j}{2})^n$.  
\end{proof}

In analogy to the billiard results stated in the introduction 
we now study the dynamics of $\mu$.  On the
space $\hC$ the dynamics are not very interesting: 
$\mu^m(Q) \to (0,\dots,0)$ as $m \to \infty$ for all $Q \in \hC$.
To get a somewhat more interesting behavior notice that
$\star$-similarity is an equivalence relation. Let $[Q] := \{ \ell Q:
\ \ell \in \C \setminus \{0\}\}$ denote the
equivalence class of $Q$.
By identifying $\star$-similar polygons we obtain the quotient space
$\hat{\hC} := \{[Q]: \ Q \in \hC \}$. 
If $Q_1$ and $Q_2$ are $\star$-similar then
so are $\mu(Q_1)$ and $\mu(Q_2)$, and thus the map $\mu$ defines a
map $\hat{\mu}$ of $\hat{\hC}$ to itself. 

We remark that the equivalence
class $[Q]$ of a polygon $Q = \sum a_i X_i$
can be represented by the vector $(a_1,\dots,a_{n-1}) \in \C^{n-1} \setminus
\{(0,\dots,0)\}$ with the identification $(a_1,\dots,a_{n-1}) \equiv
(\ell a_1,\dots,\ell a_{n-1})$ with $\ell \in \C \setminus \{0\}$.
Thus the quotient space $\hat{\hC}$ is naturally identified with the complex
projective space $\CP$.  

Let $\S$ be the unit sphere in $\hC$, i.e., 
the set of all $Q = \sum a_i X_i \in \hC$ such that $\sum |a_i|^2 = 1$.
Each equivalence class $[Q]$ intersects the unit sphere $\S$ in a
circle, i.e., $[Q] \cap \S = \{\ell Q: |\ell| = 1/\sum |a_i|^2\}$.
We define ${{\rm{dist}}}([P],[Q]) = \inf \{ d(P_0,Q_0): \ P_0 \in [P] \cap \S,
Q_0 \in [Q] \cap \S\}$, where
$d((a_1,\dots,a_{n-1}),(b_1,\dots,b_{n-1})) = (\sum |a_i -
b_i|^2)^{1/2}$ is the Euclidean
distance on $S$.   It is a simple exercise to verify that 
this defines a metric on $\hat{\hC}$.

A $\hat{\mu}$-invariant set $A \subset \hat{\hC}$ is 
%a {\em global attractor} for $\hat{\mu}$
%if, for any polygon $Q$, $dist(\hat{\mu}^n(Q), A) \to 0$ as $n \to \infty$ .
%We call $A$ 
an {\em  exponential attractor} with basin $B$ for $\hat{\mu}$  
if there exists $c,\gamma >
0$ such that $\rm{dist}(\hat{\mu}^m([Q]),A) \le c \exp(-\gamma m)$ for
all $m \ge 0$ and all $[Q] \in B$. If $B = \hat{\hC}$ then we say that
$A$ is a {\em global} exponential attractor. 

For each $j \in \{1,\dots,\lceil n/2 \rceil -1\}$ 
let $B_j \subset \hC$ be the subspace generated by
$\{X_j,X_{j+1},\dots,X_{n-j}\}$. 
Let  $B_{\lceil n/2 \rceil} := \emptyset$, and 
${\rm{Aff}} := \cup_{j=1}^{\lceil n/2 \rceil - 1} {A_j}$.
For any subset $D \subset \hC$ let $\hat{D} := \{ [P] : P \in D\}$. 
We will only use this notation for sets $D$ which are maximal in the
sense that if $\hat{D} = \hat{E}$ then $E \subset D$.

\begin{theorem}\label{nn}
For each $j \in \{1,\dots,\lceil n/2 \rceil -1\}$
the set $\hat{A_j}$ is an exponential attractor with basin 
$\hat{B_j}~\setminus~\hat{B_{j+1}}$ and thus 
%the set of affinely regular equivalence classes 
$\hat{\rm{Aff}}$ is a 
global exponential attractor for $\hat{\mu}$.
The map $\hat{\mu}$ is $n$-periodic (i.e., $\hat{\mu}^n = Id$) 
on each $\hat{A_j}$, and thus on $\hat{\rm{Aff}}$.
\end{theorem}

\begin{proof}
Suppose $Q \in B_j \setminus B_{j+1}$.  Equation (\ref{aha}) 
implies that the polygon $\mu^m(Q) = \sum_{i=1}^{n-1} (\frac{1+q^i}{2})^m a_i X_i 
= \sum_{i=j}^{n-j} (\frac{1+q^i}{2})^m a_i X_i$
is $\star$-similar to 
$$
\left (\frac{2}{1+q^j} \right )^m \mu^m(Q) = a_j X_j + \sum_{i=j+1}^{n-j-1}
\left (\frac{1+q^i}{1+q^j} \right )^m  a_i X_i 
+ \left ( \frac{1+q^{n-j}}{1+q^j} \right )^m  a_{n-j} X_{n-j}.
$$ 
Since $|(1+q^i)/(1+q^j)| < 1$ the terms in the middle sum are
exponentially small. Thus $\hat{\mu}^m ([Q])$ is exponentially close to
the equivalence class $[P_m] \in \hat{A_j}$ of the polygon 
$P_m := a_j X_j +  \left [ (1+q^{n-j})/(1+q^j) \right ]^m a_{n-j} X_{n-j}$. 
The fact that $\hat{\rm{Aff}}$ is a global attractor follows 
from this since $\hat{\hC}$ is the
union of the $\hat{B_j}$.
Lemma \ref{n} immediately implies that $\hat{\mu}$ is periodic on each 
$\hat{A_j}$.
\end{proof}

%Similarly, for each $j=2,\dots,\lceil n/2 \rceil -1$
%the equivalence class $[P]$ of affinely regular polygons $P \in A_j$
%is an exponential attractor for the restriction of 
%the map $\hat{\mu}$ to $[P]$.  
%The basin of attraction of this restriction is $[B_j] \setminus [B_{j+1}]$. 

\paragraph{Acknowledgments.}
Many thanks to Sergei Tabachnikov and the  two anonymous referees
for helpful remarks.

\bigskip
\noindent
\textbf{Serge E.  Troubetzkoy} received his B.A.\ from Yale University in
1982 and his Ph.D.\ from Stanford University in 1987. After a Brownian
career path: Leningrad, Toronto, Warwick, Bielefeld, Stony Brook,
and UAB, he settled down in Marseille where, due to  the complicated
French system, he has four affiliations (but only one salary):
Centre de Physique Th\'eorique,
F\'ed\'eration de Recherches des Unit\'es de Math\'ematique de Marseille,
Institut de Math\'ematiques de Luminy, and 
Universit\'e de la M\'editerran\'ee.  His area of research is
dynamical systems with emphasis on billiards.  He has written 
articles on four types of billiards: 
polygonal billiards, hyperbolic billiards, elliptic
billiards, and now on dual billiards.

\noindent
\textit{Centre de Physique Th\'eorique, Luminy, Case 907, F-13288 Marseille
  Cedex 9, France\\
troubetz@iml.univ-mrs.fr}

\end{document}